\documentclass[journal ]{new-aiaa}
\usepackage[utf8]{inputenc}
\usepackage{textcomp}

\usepackage{graphicx}
\usepackage{amsmath}
\usepackage[version=4]{mhchem}
\usepackage{siunitx}
\usepackage{longtable,tabularx}
\usepackage{subfig,float}
\usepackage{algorithm}
\usepackage[noend]{algpseudocode}

\usepackage[colorinlistoftodos]{todonotes}

\definecolor{orange}{rgb}{1,0.5,0}
\definecolor{amethyst}{rgb}{0.6,0.4,0.8}
\definecolor{aureolin}{rgb}{0.99,0.93,0.0}
\definecolor{awesome}{rgb}{1.0,0.13,0.32}
\definecolor{ao-green}{rgb}{0.0, 0.5, 0.0}

\newcommand{\EE}[1]{\mathbb{E}\left[#1\right]}
\newcommand{\Var}[1]{\mathbb{V}\left[#1\right]}
\newcommand{\CV}[1]{\mathbb{CV}\left[#1\right]}

\newcommand{\R}{\mathbb{R}}

\title{Rapid Aerodynamic Shape Optimization Under Parametric and Turbulence Model Uncertainty: A Stochastic Gradient Approach}

\author{Llu\'is Jofre\footnote{Assistant Professor. E-mail: lluis.jofre@upc.edu}}
\affil{Dept. Fluid Mechanics, Technical University of Catalonia - BarcelonaTech, Barcelona 08019, Spain}
\author{Alireza Doostan\footnote{Associate Professor. E-mail: alireza.doostan@colorado.edu}}
\affil{Smead Aerospace Engineering Sciences, University of Colorado, Boulder, CO 80309, USA}

\begin{document}

\maketitle

\begin{abstract}
Aerodynamic optimization is ubiquitous in the design of most engineering systems interacting with fluids. A common approach is to optimize a performance function -- subject to some constraints -- defined by a choice of an aerodynamic model, e.g., turbulence RANS model, and at nominal operating conditions. Practical experience indicates that such a deterministic, i.e., single-point, approach may result in considerably sub-optimal designs when the adopted aerodynamic model does not lead to accurate flow predictions or when the actual operating conditions differ from those considered in the design. One approach to address this shortcoming is to consider an average or robust design, wherein the statistical moments of the performance function, given the uncertainty in the operating conditions and the aerodynamic model, is optimized. However, when the number of uncertain inputs is large or the performance function exhibits significant variability, an accurate evaluation of these moments may require a large number of forward and/or adjoint solves, at each iteration of a gradient-based scheme. This, in turn, renders the design of complex aerodynamic systems computationally expensive, if not infeasible. To tackle this difficulty, we consider a variant of the stochastic gradient descent method where, in each optimization iteration, a stochastic approximation of the objective, constraints, and their gradients are generated. This is done via a small number of forward/adjoint solves corresponding to random selections of the uncertain parameters and aerodynamic model. The methodology is applied to the robust optimization of the standard NACA-0012 airfoil in a low-Mach-number turbulent flow regime and subject to parametric and turbulence model uncertainty. With a cost that is a small factor larger than that of the deterministic approach, the stochastic gradient approach significantly improves the performance (mean and variance) of the aerodynamic design for a wide range of operating conditions and turbulence models.
\end{abstract}

\section*{Nomenclature}

{\renewcommand\arraystretch{1.0}
\noindent\begin{longtable*}{@{}l @{\quad=\quad} l@{}}
$c$& airfoil chord line \\
$\textbf{C}$& stochastic constraint measures \\
$C_D$& drag coefficient \\
$C_L$& lift coefficient \\
$C_L^\star$& target lift coefficient \\
$f$& deterministic objective measure \\
$\textbf{g}$& deterministic constraint measures \\
$\textbf{h}$& stochastic optimization gradient \\
$n$& number of vector elements \\
$N$& number of optimization iterations \\
$R$& stochastic optimization objective \\
$U$& time-averaged velocity \\
$U_{\infty}$& free-stream velocity \\
$Re_{c}$& chord Reynolds number \\
$\alpha$& angle of attack \\
$\epsilon$& small user-specified number \\
$\eta$& step size / learning rate \\
$\boldsymbol{\theta}$& optimization variables \\
$\boldsymbol{\kappa}$& user-specified parameters \\
$\lambda$& robust optimization parameter \\
$\nu$& fluid kinematic viscosity \\
$\boldsymbol{\xi}$& uncertain parameters
\end{longtable*}}


\section{Introduction}  \label{sec:intro}

\lettrine{T}{he} design process of an engineering structure or device requires an appropriate selection of values for the design parameters, such that the desired performance of the system is optimized given some constraints. The specification of one or more design operating conditions allows practitioners to use deterministic optimization approaches.
For instance, in airfoil shape design, the target lift coefficient and structural properties are specified and the objective is to minimize drag under some constraints. However, the use of deterministic, single-point (DSP) formulations may lead to under-performing designs when operating away from the specific, often nominal, design conditions selected for the optimization~\cite{Huyse2002-A}. Similarly, small manufacturing imperfections or fluctuations in the flight conditions may lead to considerable changes in the airfoil performance. In practice,  the variability in the operating conditions, e.g., Reynolds or Mach number, cannot be completely eliminated. Tightening the manufacturing tolerances may prove prohibitively expensive or practically impossible to achieve; it is typically expensive to produce a precise design and impossible to maintain a pristine shape during accumulated routine flights. Beyond the aforementioned uncertainties, the optimal design depends on the choice of the aerodynamic model, e.g., a RANS model, and consequently an inaccurate model may lead to an overall sub-optimal design. 

The sensitivity of the design performance to such uncertainties provides an incentive to pursue the so-called {\it robust} designs, wherein the optimization formulation incorporates the effects of system uncertainties. Such effects may be considered by adding constraints on the failure probability of a system. This approach is known as reliability-based optimization (RBO), and mostly uses first- or second-order Taylor series expansion to approximate the moments of the limit state function with respect to the uncertain parameters~\cite{Haldar2000-B}.
Other strategies to solve RBO problems have been proposed, including polynomial chaos expansion~\cite{Eldred2011-A} and Karhunen-Lo\`eve expansion~\cite{Guest2008-A}.
On the other hand, robust optimization considers the impact of uncertainties directly on the objective and constrains through their second-order statistics~\cite{Padula2006-A}, i.e., variance or standard deviation.
Generally, the statistical moments and their gradients are estimated by means of Monte Carlo simulation, which requires many sample evaluations when the variance of the objective/constraints is large, resulting in high computational costs. On the contrary, approaches based on polynomial chaos expansions~\cite{keshavarzzadeh2016gradient} and stochastic collocation\cite{lazarov2012topology}, while effective, are limited to relatively low-dimensional uncertain inputs~\cite{doostan2009least,doostan2011non}. 

To alleviate the computational burden, and motivated by recent advances in machine learning~\cite{Bottou2018-A}, De et al.~\cite{De2020-A,de2020bi,de2021reliability} proposed stochastic gradient descent (SGD) methods~\cite{Bottou2018-A} for topology optimization of structural systems. The approach generates an unbiased, stochastic approximation of the gradients at each optimization iteration in a manner similar to the standard Monte Carlo method but with a small, e.g., $\mathcal{O}(1)$, number of realizations of the gradients (and the objective/constraints). These estimates are performed statistically independently across the optimization iterations. Through their application to several problems involving high-dimensional uncertain inputs it was demonstrated that SGD methods produce robust designs that achieve objectives similar to large sample size Monte Carlo simulation but at a cost that is only a small factor (e.g., 4) larger than solving the deterministic formulation of the optimization per iteration. 

Building upon the work in~\cite{De2020-A,de2020bi,de2021reliability}, the objectives of this work are to (i) extend the SGD approach to aerodynamic shape optimization under uncertainty, and (ii) demonstrate its efficacy in designing an airfoil subject to combined parametric and model uncertainty. In particular, we consider both average and robust design of a standard NACA-0012 airfoil in a low-Mach-number turbulent flow regime. The variability of the cruise conditions is accounted for by considering a range of Reynolds numbers. Further, to achieve a design that is less sensitive to the choice of the turbulence model, i.e., less biased by the closure assumptions/formulations, five different RANS models are considered via a discrete random variable that represent model uncertainty. The numerical results suggest the SGD approach is able to achieve designs that exhibit smaller variance in the presence of such uncertainties and with a cost that is only four times larger than that of a deterministic, single-point design.  
The rest of this manuscript is organized as follows. Section~\ref{sec:aerodynamic_shape_optimization} motivates and formally introduces the problem of aerodynamic shape optimization, together with the corresponding sources of uncertainty.
The description of the SGD strategy utilized in this work is discussed next in Section~\ref{sec:SGD}. In Section~\ref{sec:results}, the application of the SGD approach to robust optimization of the NACA-0012 airfoil is described and the corresponding results are discussed.
Finally, in Section~\ref{sec:conclusions}, the work is concluded and future directions are proposed.

\section{Aerodynamic Shape Optimization}    \label{sec:aerodynamic_shape_optimization}

This section presents the deterministic airfoil shape optimization setup with the corresponding mathematical notation and provides a discussion on the sources of uncertainty impacting the flow predictions. Subsequently, to account for these uncertainties, the robust optimization formulation considered here is introduced.   

\subsection{Deterministic Airfoil Optimization}   \label{sec:deterministic_airfoil_optimization}

The problem of optimizing the shape of an airfoil is generally formulated in terms of reducing the drag under the constraints of a minimum lift and a fixed total volume (or area in 2D) $\Omega$. The drag and lift are customarily represented by their corresponding coefficients $C_D$ and $C_L$.
The minimum required lift is also expressed by means of a target lift coefficient denoted $C_L^\star$.
Typically, the set of variables $\boldsymbol{\theta}\in \boldsymbol{\Theta}$ used to optimize the performance of the airfoil correspond to the angle of attack (AoA) $\alpha$, and the wall-normal positions of the control nodes defining the airfoil geometry.
Finally, defining the objective function $f\left( \boldsymbol{\theta}\right)$ as  $C_D$, the deterministic optimization problem is formulated as
\begin{equation}
\begin{aligned}
& \underset{\boldsymbol{\theta} \in\boldsymbol{\Theta}}{\text{minimize}}
& & C_D\left(\boldsymbol{\theta}\right)\\
& \text{subject to}
& & C_L\left(\boldsymbol{\theta}\right) \geq C_L^\star, \\
&&& \Omega = \Omega_0,    \label{eq:constrained_problem}
\end{aligned}
\end{equation}
where $\Omega_0$ is the initial volume of the airfoil.

\subsection{Sources of Uncertainty}   \label{sec:uncertainty_sources}

The aerodynamic optimization of airfoils may be affected by various sources of uncertainty.
These include aleatoric uncertainty inherent to the variability in flight conditions and model-form uncertainty as a result of the computational approaches required for rapid flow prediction in high-Reynolds-number regimes and/or complex geometrical configurations.

In the case of low-Mach-number turbulent flows, the variability in flight conditions, e.g., fluctuations in velocity and/or fluid properties, can be grouped into the chord Reynolds number $Re_{c} = U_{\infty} c / \nu$ via dimensional analysis~\cite{Jofre2020a-A}, where $U_{\infty}$ is the free-stream velocity, $c$ is the length of the airfoil's chord line, and $\nu$ is the fluid's kinematic viscosity. While this uncertainty can be characterized based on experimental data obtained from flight tests and databases, it cannot be reduced as it depends on the properties of the ambient fluid in which the airfoil is operating; these are generally dictated by large-scale atmospheric phenomena and local weather conditions.

At present, flow prediction over aerodynamic airfoils under turbulent conditions utilizing scale-resolving computational approaches, such as direct numerical simulation (DNS) and large-eddy simulation (LES), is expensive. In addition to the forward flow solves, gradient-based, shape optimization techniques require information about the sensitivity of the objective function and constraints with respect to the design variables computed, e.g., backward/adjoint computations. Therefore, these high-fidelity computational approaches are not a practical option in  optimization of large-scale applications. As such, design strategies often rely on computationally cheaper solutions based on Reynolds-averaged Navier-Stokes (RANS)~\cite{Wilcox1998-B} and detached-eddy simulation (DES)~\cite{Spalart2009-A}. The study and optimization of complex turbulent flows by means of RANS simulations have become routine in many scientific and engineering applications.
The underlying time-averaging operation of the approach enables significant reduction of the computational requirements by resolving only the integral flow scales dictated by the boundary conditions of the problem.
However, the Reynolds stresses and their effects on the resolved flow field are not negligible, and therefore require additional modeling.
Consequently, the assumptions made in the closure formulations become important sources of model-form uncertainty that impact the QoIs~\cite{Emory2013-A,Jofre2018-A,Jofre2019-A}.
As a result, many RANS models tailored (initially) for specific flow regimes and configurations exist in the literature~\cite{Wilcox1998-B}.
In the present study, therefore, this type of uncertainty is also introduced into the optimization problem by considering fundamentally different RANS models which are randomly selected based on a discrete random variable with with equal probability masses.

It is important to note that, over the past years, multifidelity (MF) strategies to reduce the cost of the outer-loop problem by combining the accuracy of high-fidelity (HF) models, e.g., DNS and LES, with the speedup achieved by low-fidelity (LF) representations, e.g., RANS, DES and analytical functions/correlations, have been extensively developed~\cite{Giles2008-A,Fairbanks2017-A,Peherstorfer2018-A,Hampton2018-A,Adcock2020-A} and applied to propagate uncertainty in large-scale multiphysics turbulent flows~\cite{Jofre2017-A,Fairbanks2020-A,Jofre2020b-A,Valero2021-A} and optimization problems~\cite{Leary2003-A,Robinson2006-A,Forrester2007-A,Lam2015-A}. As demonstrated in \cite{de2020bi}, the SGD approach of this study may be extended to incorporate MF aerodynamic models, a direction we leave for a future study. 


\subsection{Robust Optimization}   \label{sec:robust_optimization}

Developing optimization methods that result in more robust designs is of significant interest in engineering practice.
The term robustness refers to a variety of goals: (i) identify designs that minimize the variability of a manufactured product in the presence of uncertainty, (ii) mitigate the detrimental effects of the worst-case performance, or (iii) obtain a uniform improvement of the performance over the entire range of operating conditions. 

Let $\boldsymbol{\theta} \in \R^{n_{\boldsymbol{\theta}}}$ denote the vector of optimization variables and $\boldsymbol{\xi} \in \R^{n_{\boldsymbol{\xi}}}$ the vector of random variables characterizing system uncertainties.
Let $f\left( \boldsymbol{\theta} ; \boldsymbol{\xi} \right) : \R^{n_{\boldsymbol{\theta}}} \times \R^{n_{\boldsymbol{\xi}}} \rightarrow \R$ denote the performance function for an instance of $\boldsymbol{\theta}$ and a realized value of $\boldsymbol{\xi}$.
Similarly, let $\textbf{g}\left( \boldsymbol{\theta} ; \boldsymbol{\xi} \right) : \R^{n_{\boldsymbol{\theta}}} \times \R^{n_{\boldsymbol{\xi}}} \rightarrow \R^{n_{\textbf{g}}}$ be the vector of $n_{\textbf{g}}$ real-valued constraints. We say that $\left( \boldsymbol{\theta} ; \boldsymbol{\xi} \right)$ satisfies the constraints if $\textbf{g}\left( \boldsymbol{\theta} ; \boldsymbol{\xi} \right) \leq 0$, and refer to positive values of $\textbf{g}\left( \boldsymbol{\theta} ; \boldsymbol{\xi} \right)$ as constraint violations.

The robust optimization objective considered here is defined as a combination of its expected value and variance as
\begin{equation}
  R\left( \boldsymbol{\theta} \right) = \EE{ f\left( \boldsymbol{\theta} ; \boldsymbol{\xi} \right) } + \lambda_{R} \Var{ f\left( \boldsymbol{\theta} ; \boldsymbol{\xi} \right) },   \label{eq:objective}    
\end{equation}
where $\EE{\cdot}$ and $\Var{ \cdot }$ denote the mathematical expectation and variance of their arguments, respectively.
Similarly, the constraint violation is expressed by
\begin{equation}
  \textbf{C}_{j}\left( \boldsymbol{\theta} \right) = \EE{ \textbf{G}_{j}\left( \boldsymbol{\theta} ; \boldsymbol{\xi} \right) } + \lambda_{C,j} \Var{ \textbf{G}_{j}\left( \boldsymbol{\theta} ; \boldsymbol{\xi} \right) }, \quad j = 1, \dots, n_{\textbf{g}},   \label{eq:constraint}    
\end{equation}
where $\textbf{G}_{j}\left( \boldsymbol{\theta} ; \boldsymbol{\xi} \right) = \{ \textrm{max} \left[ 0, \textbf{g}_{j}\left( \boldsymbol{\theta} ; \boldsymbol{\xi} \right) \right] \}^{2}$ for $j = 1, \dots, n_{\textbf{g}}$,~\cite{Griva2009-B}.
Note that while $f$ and $\textbf{g}$ are random variables, respectively, $R$ and $\textbf{C}_j$ are scalars as they are associated with moments taken with respect to the probability measure of $\boldsymbol{\xi}$.
The parameters $\lambda_{R}, \lambda_{C,j} \geq 0$ denote, respectively, the importance of variations in $f$ or $\textbf{g}$ relative to their means. The effect of adding the variance to the objective and constraints is to aim for a design that shows relatively smaller variations in the two responses, $f$ and $\textbf{g}$, even in the presence of uncertainty.
Hence, we denote this formulation as \textit{robust optimization} formulation for $\lambda_{R}, \lambda_{C,j} > 0$, and {\it average optimization} when $\lambda_{R}, \lambda_{C,j} = 0$. Although one may use different values for $\lambda_{R}$ and $\lambda_{C,j}$ in Eqs.~(\ref{eq:objective}) and~(\ref{eq:constraint}), we restrict the analysis to utilizing the same value in the application studied in Section~\ref{sec:optimization_results}; viz. the two parameters are considered to be equal and denoted by $\lambda$.

Finally, we are interested in solving the unconstrained optimization problem expressed as
\begin{equation}
\begin{aligned}
& \underset{\boldsymbol{\theta} \in\boldsymbol{\Theta}}{\text{minimize}}
& & R\left( \boldsymbol{\theta} \right) + \boldsymbol{\kappa}^{\intercal}\textbf{C}\left( \boldsymbol{\theta} \right),    \label{eq:unsconstrained_problem}
\end{aligned}
\end{equation}
where $\boldsymbol{\kappa}$ is a user-specified vector of parameters that penalizes against constraint violation.
Generally, an extensive search for $\boldsymbol{\kappa}$ is not computationally feasible.
Hence, based on a few preliminary tests, the values were determined to $\kappa_{j} = 1$ for the application studied in this paper.

\section{Stochastic Gradient Descent Approach}    \label{sec:SGD}

In standard Monte Carlo approaches, $R\left( \boldsymbol{\theta} \right)$, $\textbf{C}\left( \boldsymbol{\theta} \right)$, and their gradients are approximated utilizing $N\gg 1$, e.g., $N = 10^3$, forward and adjoint solves of the model for specific values of the design variables $\boldsymbol{\theta}$ and $N$ realizations of $\boldsymbol{\xi}$.
However, evaluating $f\left( \boldsymbol{\theta},\boldsymbol{\xi} \right)$, $\textbf{g}\left( \boldsymbol{\theta},\boldsymbol{\xi} \right)$, and their gradients, at each iteration $N$ times may be computationally expensive for problems with many degrees of freedom.
Instead, following \cite{De2020-A,de2020bi,de2021reliability}, we use $n\sim\mathcal{O}(1)\ll N$, e.g., $n=2, 4, 8$,  random samples to estimate the expectations, variances, and gradients of the objective and constraints. Such stochastic approximations are performed independently across optimization iterations; that is, an independent set of $n$ samples of $\boldsymbol{\xi}$ is used at each iteration.

The SGD method shown in Algorithm~\ref{alg:sgd} uses a single, i.e., $n=1$, realization of $\boldsymbol{\xi}$ to update the design at the $k$th iteration as
\begin{equation}
\boldsymbol{\theta}_{k+1} = \boldsymbol{\theta}_{k} - \eta\textbf{h}_k
\end{equation}
where the gradient $\textbf{h}_k$ is defined as
\begin{equation}
  \textbf{h}_k = \nabla f\left( \boldsymbol{\theta}_k, \boldsymbol{\xi}_i \right) + \boldsymbol{\kappa}^{\intercal}\nabla\textbf{G}\left( \boldsymbol{\theta}_k, \boldsymbol{\xi}_i \right), \label{eq:sgd_gradient}
\end{equation}
and $\eta$ is the step size, also known as learning rate.
Hence, the computational cost of SGD is relatively low, specifically same as the deterministic variant of the optimization at each iteration. However, the convergence of the standard SGD with $n=1$ can be slow since it may not follow the descent direction at every iteration due to the large variance of the stochastic gradients in Eq.~(\ref{eq:sgd_gradient}).

\begin{algorithm}[htb]
\caption{Standard stochastic gradient descent (SGD) method~\cite{Bottou2018-A}}	\label{alg:sgd}
\begin{algorithmic}[1]
\State Given $\eta$
\State Initialize $\boldsymbol{\theta}_1$
\State \textbf{for} $k=1,2,\dots$ \textbf{do}
\State \quad Compute $\textbf{h}_k := \textbf{h}\left( \boldsymbol{\theta}_k \right)$
\State \quad Set $\boldsymbol{\theta}_{k+1}$ $\gets$  $\boldsymbol{\theta}_{k} - \eta\textbf{h}_k$
\State \textbf{end for}
\end{algorithmic}
\end{algorithm}

A straightforward extension of the standard SGD method, known as \textit{mini-batch gradient descent}~\cite{Bottou2018-A},  is to use $n>1$ random samples at each iteration to estimate the gradient $\textbf{h}_k$. In the past few years, several modifications to the SGD have been proposed to improve its convergence, e.g., the Adaptive Subgradient (AdaGrad)~\cite{Duchi2011-A}, Adadelta~\cite{Zeiler2012-A}, Adaptive Moment Estimation (Adam)~\cite{Kingma2014-A}, Stochastic Average Gradient (SAG)~\cite{Roux2012-A}, and Stochastic Variance Reduced Gradient (SVRG)~\cite{Johnson2013-A}. The SGD algorithm employed in this work is the AdaGrad~\cite{Duchi2011-A}, described in the next subsection and summarized in Algorithm~\ref{alg:adagrad}, where historical information about the gradients is used to modify the updates of the optimization variables.

\subsection{Adaptive Subgradient Method}    \label{sec:AdaGrad}

AdaGrad dampens the movements along directions with historically large gradients, adapting in this way the learning rate and facilitating a faster convergence.
In this algorithm, at iteration $k$, the following auxiliary variable is computed
\begin{equation}
  \textbf{a}_{k,j} = \sum_{i=1}^{k}\textbf{h}_{i,j}^{2}, \quad j = 1, \dots, n_{\boldsymbol{\theta}}, \label{eq:auxiliary_variable}
\end{equation}
which is then used in the update rule
\begin{equation}
  \boldsymbol{\theta}_{k+1} = \boldsymbol{\theta}_k - \eta\textbf{a}_{k}^{-1/2}\textbf{h}_k, \label{eq:update_rule}
\end{equation}
where the vector multiplication is performed component-wise.
To avoid division by zero, a small number $\epsilon = 10^{-8}$ is incorporated in the denominator of the update, as shown in Algorithm~\ref{alg:adagrad}.

\begin{algorithm}[htb]
\caption{Adaptive Subgradient (AdaGrad) method~\cite{Duchi2011-A}}	\label{alg:adagrad}
\begin{algorithmic}[1]
\State Given $\eta$
\State Initialize $\boldsymbol{\theta}_1$
\State Initialize $\textbf{a} := \textbf{0}$, having the same dimension as $\boldsymbol{\theta}_1$
\State \textbf{for} $k=1,2,\dots$ \textbf{do}
\State \quad Compute $\textbf{h}_k := \textbf{h}\left( \boldsymbol{\theta}_k \right)$
\State \quad Set $\textbf{a}_{j} := \textbf{a}_j + \textbf{h}_{k,j}^{2}, \quad j = 1, \dots, n_{\boldsymbol{\theta}}$
\State \quad Set $\boldsymbol{\theta}_{k+1,j}$ $\gets$  $\boldsymbol{\theta}_{k,j} - \eta\frac{\textbf{h}\left(\boldsymbol{\theta}_{k,j}\right)}{\sqrt{\textbf{a}_j}+\sqrt{\epsilon}}, \quad j = 1, \dots, n_p$
\State \textbf{end for}
\end{algorithmic}
\end{algorithm}

\section{Optimization Experiments}  \label{sec:results}

The performance of the SGD strategy introduced in Section~\ref{sec:SGD}, specifically AdaGrad, is discussed in the subsections below by considering the robust optimization of an airfoil in a low-Mach-number turbulent flow regime.
The problem and uncertainty sources are described first in Section~\ref{sec:problem}, and the results are presented and discussed next in Section~\ref{sec:optimization_results}.

\subsection{Problem Description}    \label{sec:problem}

The optimization problem selected to analyze the performance of the strategy presented in this work considers the NACA-0012 airfoil with sharp trailing edge as the baseline geometry.
The objective is to minimize the drag coefficient $C_D$ at cruise conditions under uncertainty in a low-Mach-number turbulent flow regime, while providing a lift coefficient $C_L$ above a specified value $C_L^\star = 0.375$, characteristic of typical values found for commercial transport airliners~\cite{Bertin2002-B}.

The Reynolds number of the problem is defined as $Re_{c} = U_{\infty} c / \nu$, where $U_{\infty}$ is the free-stream velocity, $c$ is the length of the airfoil's chord line, and $\nu$ is the kinematic viscosity of the fluid.
The design variables in the optimization problem are given by the vertical positions of $10\times2\times2 = 40$ equidistant free-form deformation (FFD) nodes~\cite{Sederberg1986-P} and the angle of attack $\alpha$.
Three airfoil profile points are in locked positions, one at the leading edge and a double control node at the trailing edge. In addition, the volume of the airfoil is kept fixed by imposing the constraint $\Omega = \Omega_0$ into the optimization problem.

The RANS equations for the flow are discretized and solved on a two-dimensional (2-D) O-type body-fitted mesh utilizing the OpenFOAM~\cite{OpenFOAM} computational package.
A grid refinement study based on the Spalart-Allmaras RANS model was performed to minimize the relative variation of $C_D$ and $C_L$ below $5\%$ for the baseline airfoil geometry.
The resulting final mesh utilized for the optimizations is composed of approximately 25000 grid points clustered around the airfoil with the first grid point in the wall-normal direction at $\hat{n}^+\sim 1$ (wall units), and with the far-field boundaries located at 50 chord lengths.
The sensitivities of both drag and lift coefficients with respect to the design parameters are efficiently calculated utilizing a discrete adjoint solver provided by the DAFoam~\cite{He2020-A} framework.

Two combined sources of uncertainty (i.e., random inputs) are considered.
First, an aleatoric uncertainty resulting from the inherent variability of the cruise conditions, e.g., changes in fluid density and viscosity and free-stream velocity, grouped within a Reynolds number in the range $Re_c = [10^6:10^7]$.
Second, uncertainty arising from the utilization of RANS turbulence models to close the systems of equations describing the flow field.
In this regard, five different RANS models are selected for this study: (1) Spalart-Allmaras~\cite{Spalart1992-A}, (2) $k$-$\epsilon$~\cite{Launder1974-A}, (3) realizable $k$-$\epsilon$~\cite{Shih1995-A}, (4) shear stress transport (SST) $k$-$\omega$~\cite{Menter1993-A}, and (5) Langtry-Menter SST $k$-$\omega$~\cite{Langtry2009-A}; a discrete random variable with equal probability has been utilized for sampling over the RANS model space.

\subsection{Optimization Results}   \label{sec:optimization_results}

\begin{figure}[hbt!]
\centering
\subfloat[]{\includegraphics[width=0.495\textwidth]{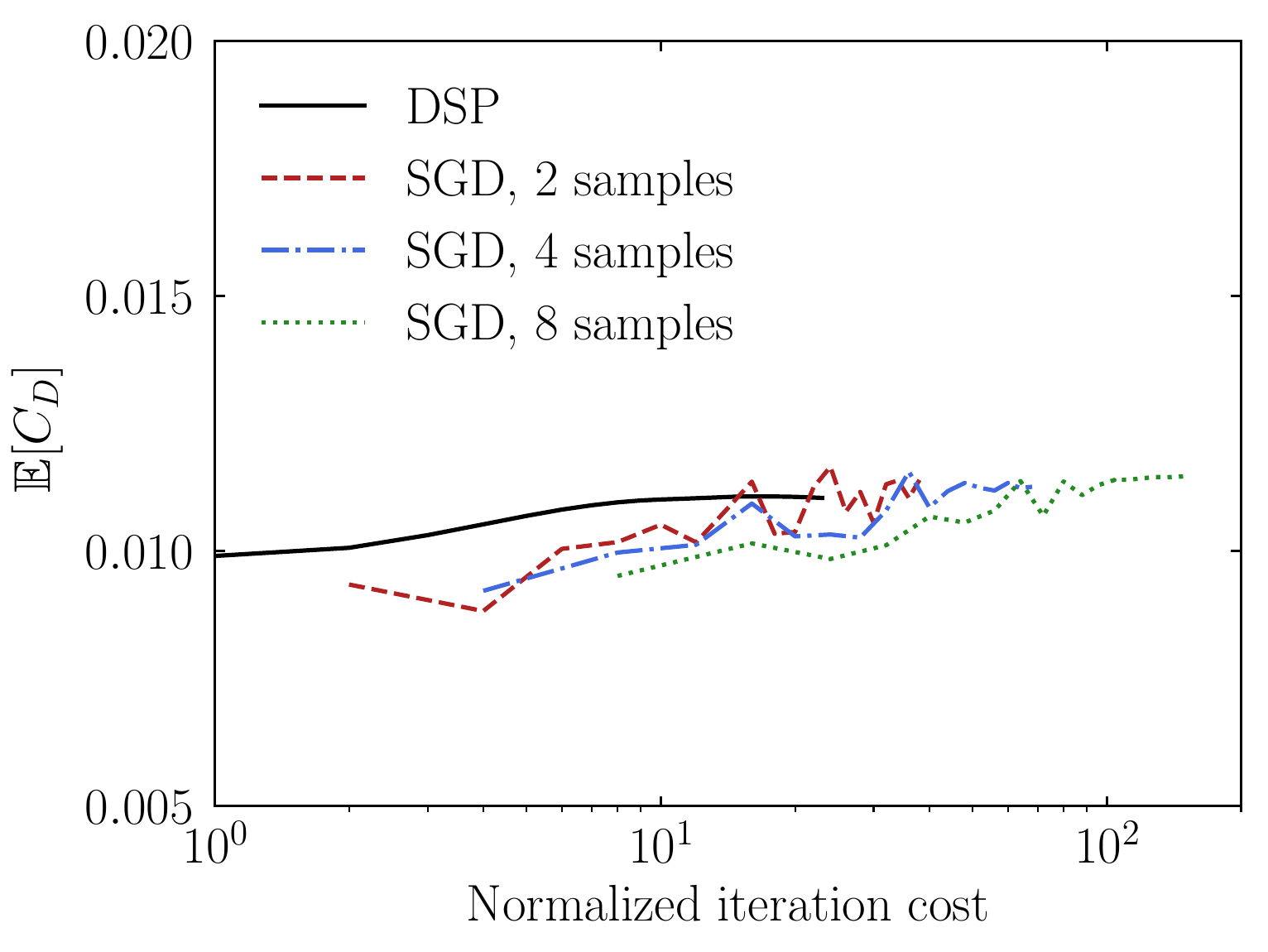}}
\hspace{3mm}
\subfloat[]{\includegraphics[width=0.479\textwidth]{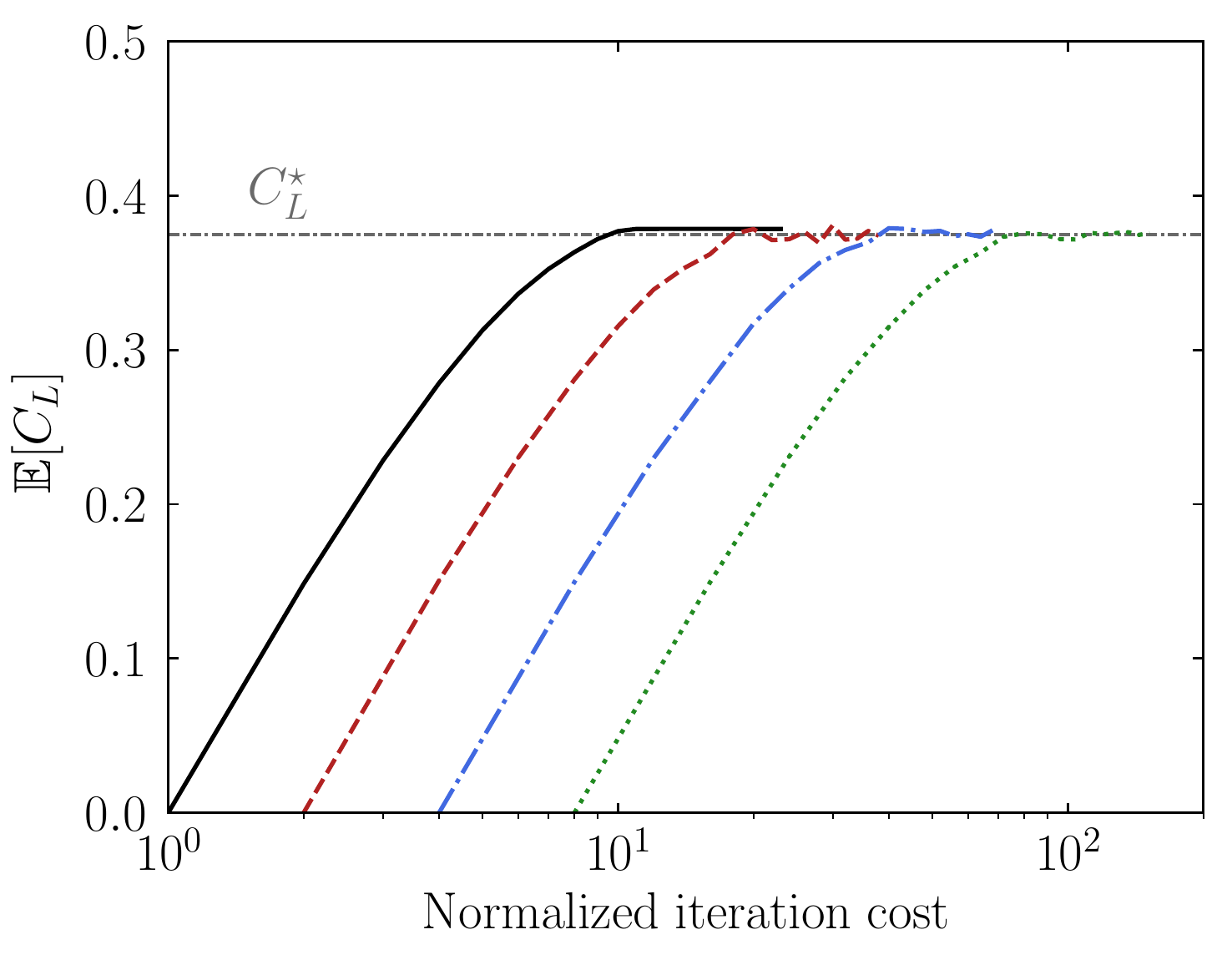}}
\caption{Evolution of mean drag (a) and lift (b) coefficients along the optimization iterations expressed in terms of normalized cost for the DSP and average (using $n=2, 4, 8$ samples per iteration) design strategies. The horizontal dot-dash line in (b) indicates the target lift coefficient $C_L^\star = 0.375$.}    \label{fig:iteration_evolution}
\end{figure}

The convergence behavior of the DSP and average optimization methodologies is studied first in Fig.~\ref{fig:iteration_evolution} by depicting the mean drag (a) and lift (b) coefficients as a function of normalized cost (or equivalent number of iterations) defined as the iteration number multiplied by the number of samples per iteration.
The DSP approach utilizes one single set of fixed values for the optimization process consisting of $Re_c = 5\cdot10^6$ and the Spalart-Allmaras RANS model (option 1), while the average methodology advances the optimization procedure by randomly sampling the stochastic inputs; viz. Reynolds numbers in the range $Re_c = [10^6:10^7]$ and RANS models from options 1, 2, 3, 4, and 5.
In particular, three cases are studied for the average design by utilizing $n\in\{2,4,8\}$ samples per optimization iteration with $\lambda = 0$.
Four main observations can be inferred from the figure.
First, the DSP optimization converges slightly faster than the average design due to their inherent differences: (i) DSP utilizes one sample per iteration, whereas average design requires more than one sample per iteration to collect statistical information and (ii) the input values are kept constant for DSP during the optimization. Second, the average design curves present oscillations, but the main trends are similar to the DSP approach.
In fact, smoother curves are obtained for average design when increasing the number of samples per iteration $n$, however, at larger computational costs. 
Third, the optimal AdaGrad performance for this problem is to utilize 4 samples per iteration: (i) 2 samples per iteration is the faster option, but it is not stable enough; and (ii) 8 samples per iteration provides slightly smoother results, but it is approximately $2\times$ more expensive than utilizing 4 samples per iteration.
Fourth, as expected and corroborated by the parallel increase of $C_D$ with $C_L$ in the figures, there is a trade-off between reducing drag and providing a minimum lift value.

\begin{figure}[hbt!]
\centering
\subfloat[]{\includegraphics[width=0.475\textwidth]{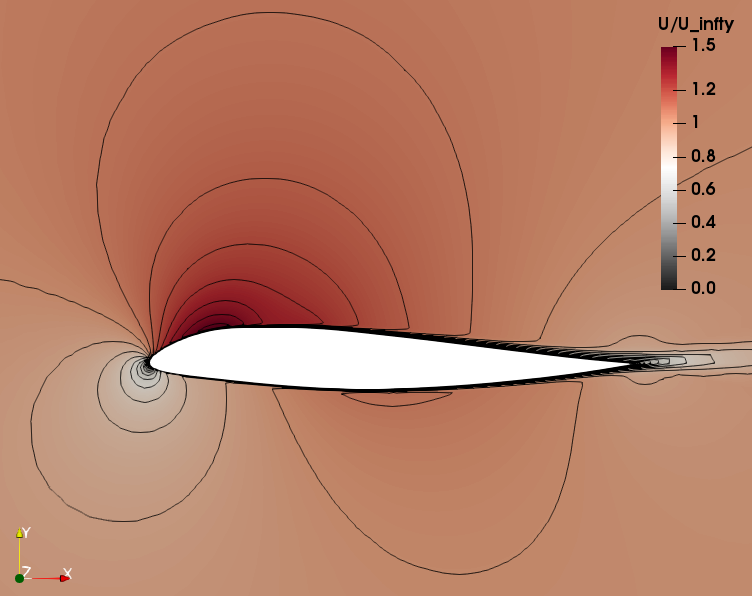}}
\hspace{3mm}
\subfloat[]{\includegraphics[width=0.475\textwidth]{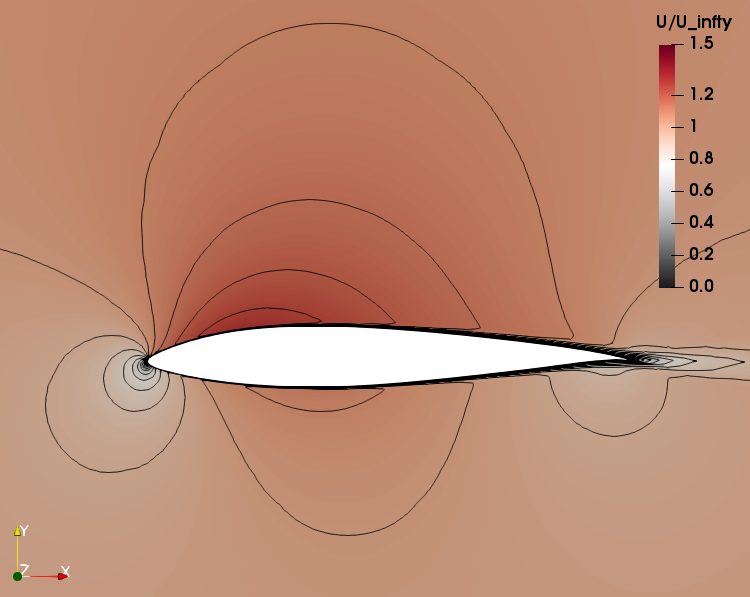}}
\vspace{3mm}
\subfloat[]{\includegraphics[width=0.475\textwidth]{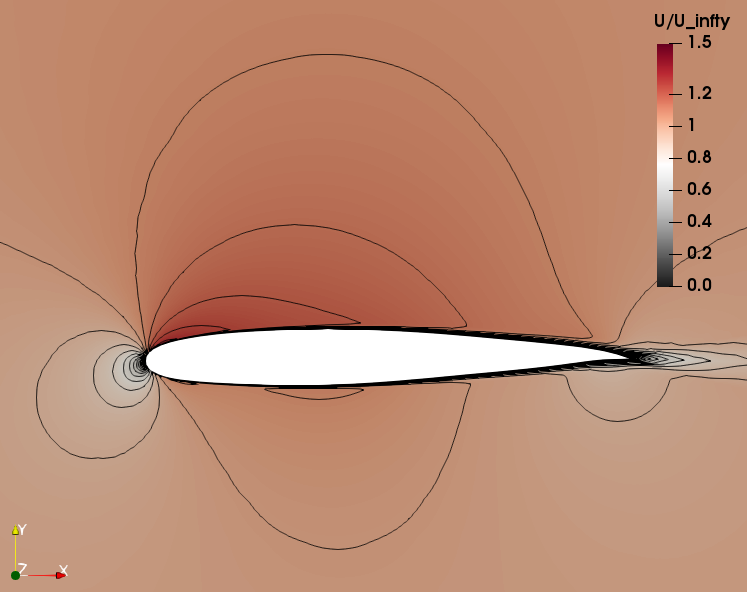}}
\hspace{3mm}
\subfloat[]{\includegraphics[width=0.475\textwidth]{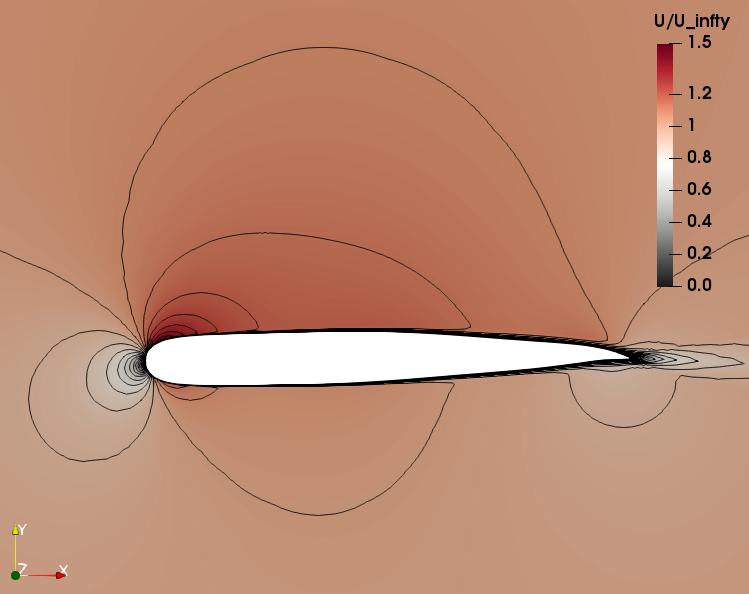}}
\caption{Normalized time-averaged velocity fields, i.e., $U/U_\infty$, obtained from the optimized airfoils using the DSP (a) and average/robust [with $\lambda = 0$ (b), $\lambda = 10$ (c), $\lambda = 100$ (d)] design strategies.}    \label{fig:optimized_airfoils}
\end{figure}

The shapes of the optimized airfoils obtained from utilizing the DSP, average, and robust designs are shown in Fig.~\ref{fig:optimized_airfoils}, together with the resulting time-averaged velocity distributions normalized by the free-stream velocity.
For the average and robust designs, 4 samples per iteration have been considered with three different values for $\lambda = 0, 10, 100$.
Following the notation in Fig.~\ref{fig:optimized_airfoils}, the corresponding angles of attack are: (a) $\alpha = 2.03^\circ$, (b) $\alpha = 1.50^\circ$, (c) $\alpha = 1.57^\circ$, and (d) $\alpha = 1.65^\circ$.
As depicted in the figure, the DSP strategy generates a significantly more asymmetric shape, with a finer leading edge, than the average- and robust-optimized airfoils.
In particular, while the region of the trailing edge does not vary much between airfoils, the leading edge gets consistently blunter with increasing $\lambda$ to increase the robustness of the final design.
The angle of attack is largest for the DSP optimization and increases with $\lambda$.
However, the difference between the $\alpha$ values obtained are not significant (in the order of half a degree).
The reason is that $C_L$ is largely sensitive to $\alpha$ in comparison to the shape of the airfoil.
Therefore, for all optimization cases the angle of attack increases until the target $C_L^\star$ value is obtained, the point at which the optimization focus is then shifted toward predominantly reducing $C_D$.
In terms of normalized time-averaged velocity fields, the DSP optimization generates a solution with a concentrated region of large velocities at the suction side of the leading edge and extending a substantial distance in the $y$-direction.
Consequently, the lift force is highly concentrated at the front part of the airfoil.
On the contrary, the average and robust strategies optimize the airfoil to generate a region of larger relative velocities that is more distributed along the entire suction side, especially for large $\lambda$ values.

\begin{figure}[hbt!]
\centering
\subfloat[]{\includegraphics[width=0.49\textwidth]{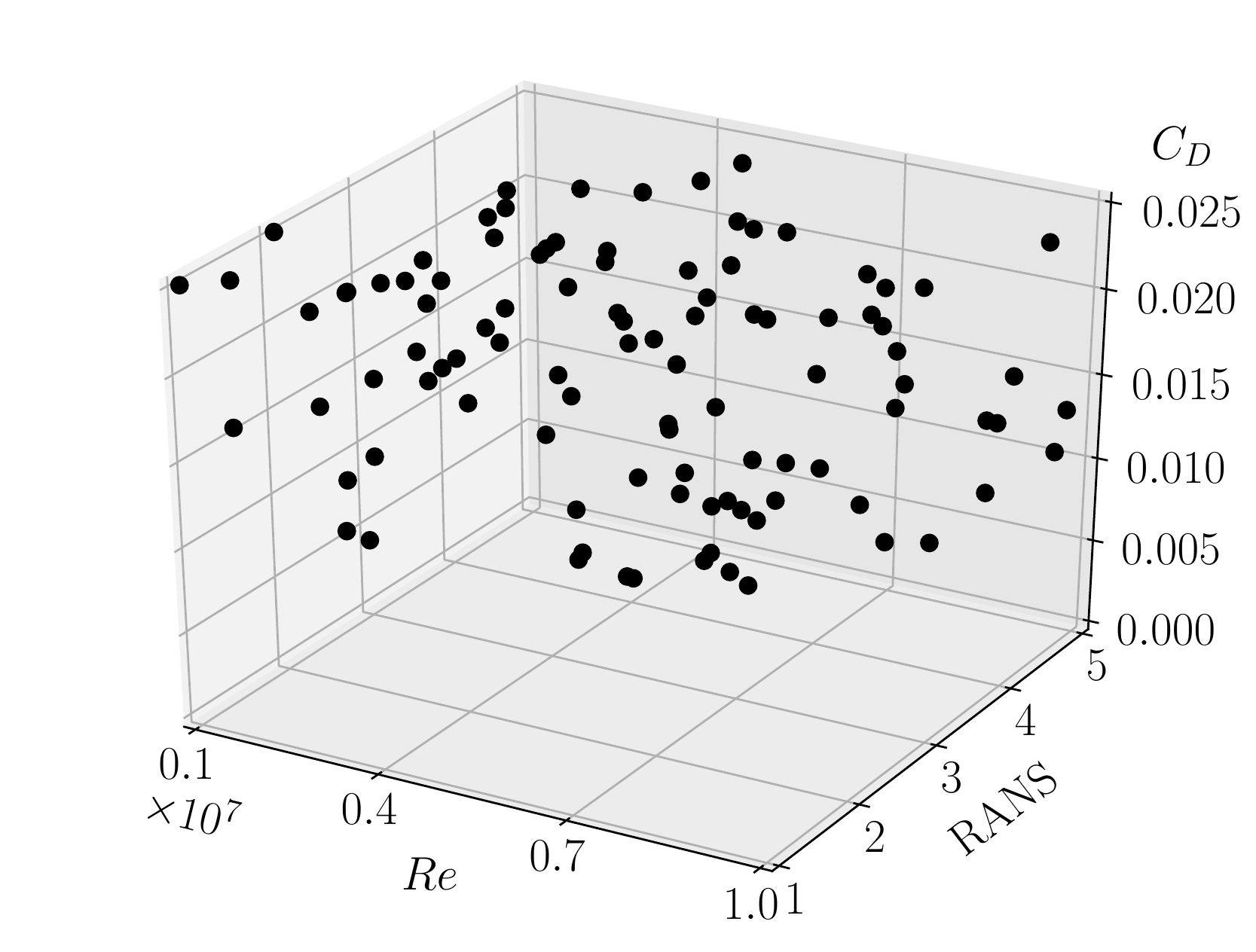}}
\subfloat[]{\includegraphics[width=0.49\textwidth]{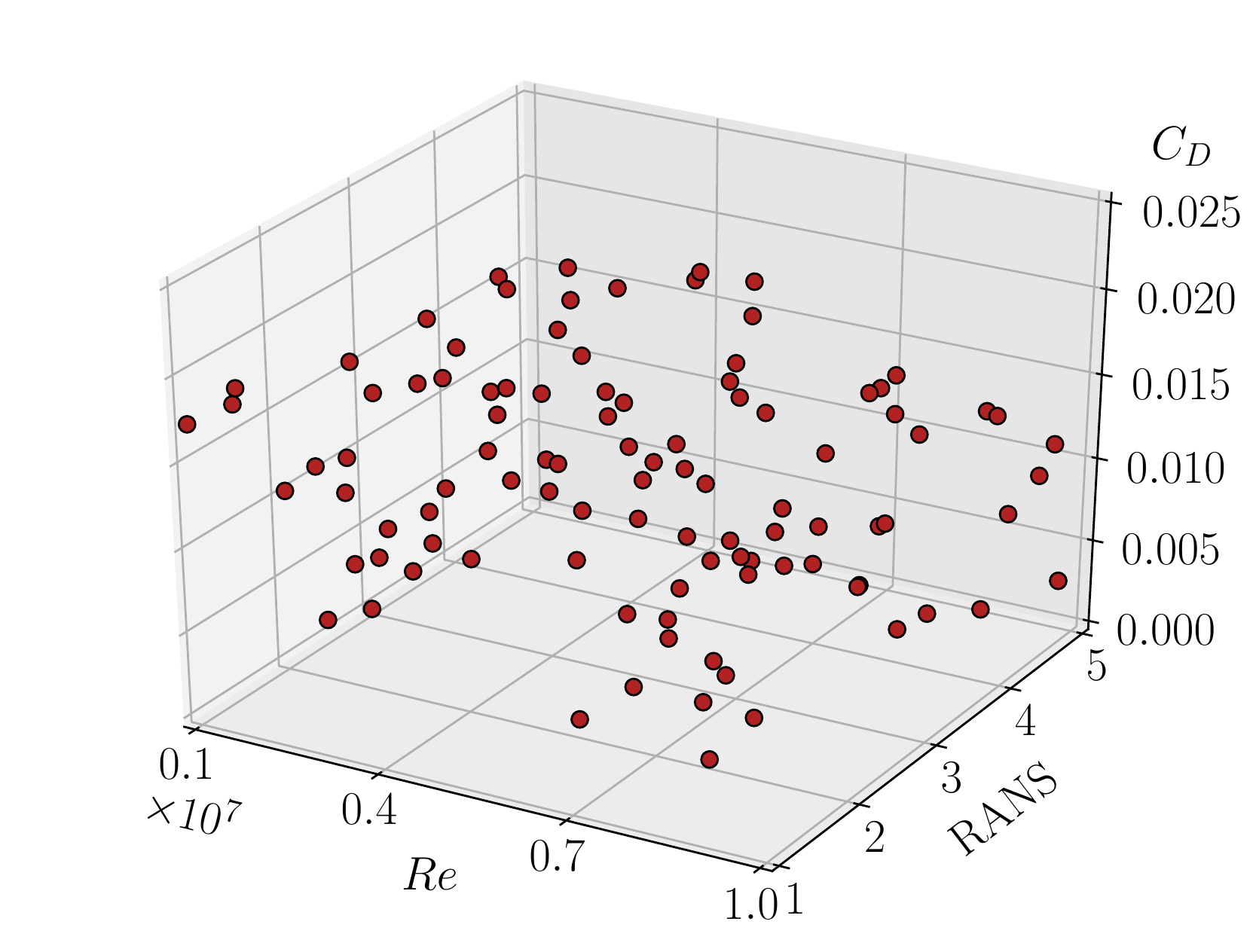}}
\\
\subfloat[]{\includegraphics[width=0.49\textwidth]{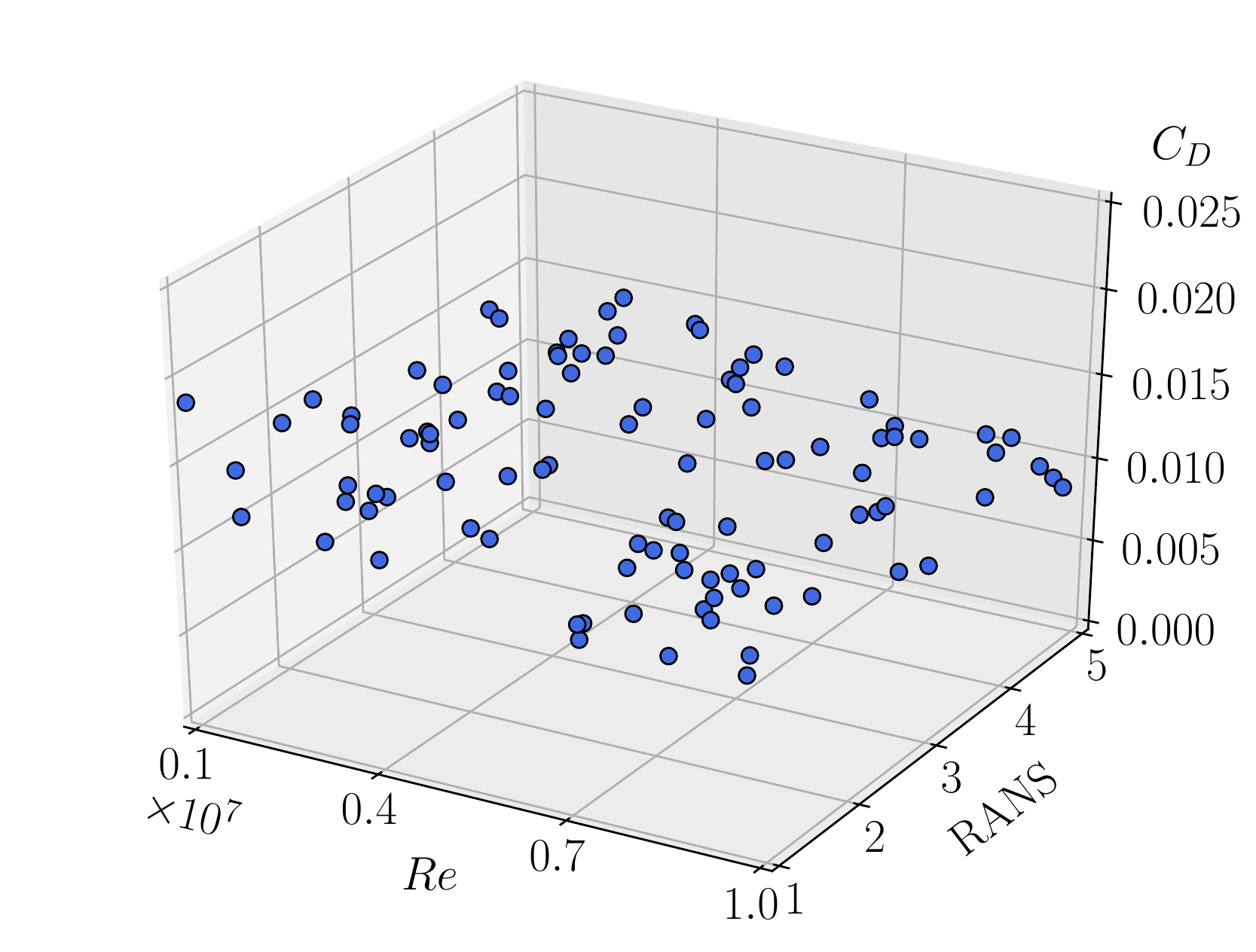}}
\subfloat[]{\includegraphics[width=0.49\textwidth]{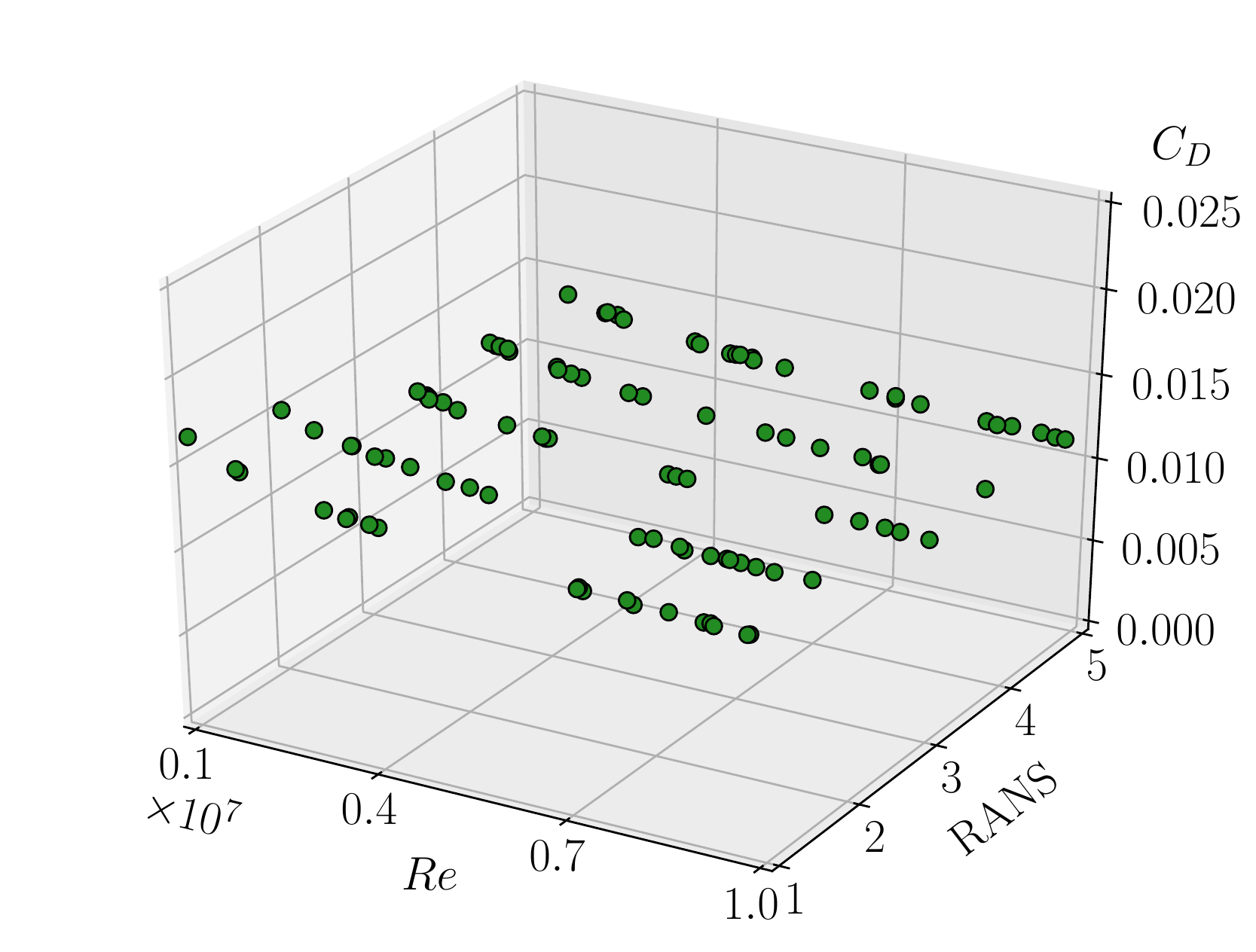}}
\caption{Drag coefficient distributions for the parameter space study obtained from the optimized airfoils using the DSP (a) and average/robust [with $\lambda = 0$ (b), $\lambda = 10$ (c), $\lambda = 100$ (d)] design strategies.}    \label{fig:parameter_space_C_D}
\end{figure}

The performance based on the drag and lift coefficients of the optimized designs described above is assessed by sampling the input parameter space by computing the same hundred randomly generated samples for each airfoil.
The resulting data is depicted in Fig.~\ref{fig:parameter_space_C_D} and summarized in Table~\ref{tab:performance_summary} in terms of mean $\EE{\cdot}$, variance $\Var{ \cdot }$, and coefficient of variation $\mathbb{CV}\left[ \cdot \right]$.
The discussion of the results can be organized in 4 major points.
First, the DSP approach provides the smallest $C_D = 9.6\cdot10^{-3}$ when optimizing for one set of $Re_c$ value and RANS model, but provides a $\EE{C_D} = 1.8\cdot10^{-2}$ for the parameter space study with $\CV{C_D} = 25\%$.
In addition, $\EE{C_L} = 4.0\cdot10^{-1}$ is above the threshold value $C_L^\star = 0.375$, but not the entire stochastic distribution since $\EE{C_L} - \Var{C_L}^{1/2} = 0.346 < C_L^\star$.
Second, the average design ($\lambda = 0$) is able to significantly reduce the mean drag coefficient to $\EE{C_D} = 1.1\cdot10^{-2}$, which is a reduction of $1.6\times$, with a larger coefficient of variation of $\CV{C_D} = 39\%$ due to a reduction of the mean value with respect to the DSP approach.
Moreover, this strategy is able to maintain $\EE{C_L} - \Var{C_L}^{1/2} = 0.380 > C_L^\star$ by increasing the mean lift coefficient value to $\EE{C_L} = 4.2\cdot10^{-1}$.
Third, the robust design with $\lambda = 10$ slightly increases the mean drag coefficient value to $\EE{C_D} = 1.2\cdot10^{-2}$, while largely reducing the coefficient of variation to $\CV{C_D} = 17\%$ (a reduction of $2.3\times$ with respect to the average design).
Additionally, the value of $\EE{C_L}$ reduces together with its $\Var{C_L}$, maintaining the requirement to $\EE{C_L} - \Var{C_L}^{1/2} = 0.376 \approx C_L^\star$.
Fourth, the robust design with $\lambda = 100$ keeps slightly increasing the mean drag coefficient value, but is able to significantly reduce the coefficient of variation to $\CV{C_D} = 8\%$ by largely decreasing $\CV{C_D}$, as clearly observed by the small variation of the data in Fig.~\ref{fig:parameter_space_C_D}(d).

\begin{table}
 \caption{Perfomance of the optimized airfoils obtained by the DSP and average/robust design strategies.} \label{tab:performance_summary}
 \centering
 \begin{tabular}{lcccccc}
 \hline
 Strategy & $\EE{C_D}$ & $\Var{C_D}$ & $\CV{C_D}$ & $\EE{C_L}$ & $\Var{C_L}$ & $\CV{C_L}$ \\\hline
 DSP                & $1.8\cdot10^{-2}$ & $2.0\cdot10^{-5}$ & $2.5\cdot10^{-1}$ & $4.0\cdot10^{-1}$ & $2.5\cdot10^{-3}$ & $1.3\cdot10^{-1}$ \\
 Average Design, $\lambda=0$   & $1.1\cdot10^{-2}$ & $1.9\cdot10^{-5}$ & $3.9\cdot10^{-1}$ & $4.2\cdot10^{-1}$ & $1.9\cdot10^{-3}$ & $1.0\cdot10^{-1}$ \\
 Robust Design, $\lambda=10$  & $1.2\cdot10^{-2}$ & $4.1\cdot10^{-6}$ & $1.7\cdot10^{-1}$ & $4.0\cdot10^{-1}$ & $4.0\cdot10^{-4}$ & $5.0\cdot10^{-2}$ \\
 Robust Design, $\lambda=100$ & $1.3\cdot10^{-2}$ & $1.1\cdot10^{-6}$ & $8.2\cdot10^{-2}$ & $3.9\cdot10^{-1}$ & $3.9\cdot10^{-4}$ & $5.0\cdot10^{-2}$ \\ 
 \hline
 \end{tabular}
\end{table}


\section{Summary \& Conclusions}   \label{sec:conclusions}

Deterministic single- and multi-point optimization strategies generally suffer from significant performance decay when aerodynamic systems operate away from the design conditions selected.
A robust solution is to perform the optimization considering the uncertainty sources by defining objectives and constraints with statistical moments of geometric and physical QoIs.
However, the necessity to compute forward and backward flow solves to accurately calculate statistical moments and their gradients pose the problem typically infeasible from a computational cost perspective.
This work, therefore, explores the acceleration strategy of generating a stochastic approximation of the objective, constraints, and their gradients by computing a small number of solves per optimization iteration.
The strategy is combined with the AdaGrad SGD optimization method, which has gained attention over the past years for data-driven scientific and engineering applications.

The methodology presented has been assessed by performing a robust optimization of the NACA-0012 airfoil in a low-Mach-number turbulent flow regime.
The results obtained indicate that, in comparison to a deterministic, single-point approach, the strategy analyzed is able to reduce the mean drag coefficient by approximately $2\times$ for a relatively wide range of operating conditions, while satisfying the constraints of the problem in terms of target lift coefficient and fixed airfoil volume.
In addition, the strategy presented enables to straightforwardly control, utilizing an external tunable parameter, the degree of robustness incorporated into the aerodynamic design as demonstrated by the reduction of variability in the response surface of the parameter space study.
Consequently, the example analyzed shows that the use of stochastic gradients estimated with a small number of random samples at each iteration (4 samples in the case studied) produces meaningful designs in a computationally affordable manner.

Future work will encompass the application of the strategy studied to the robust optimization of aerodynamic problems involving other flow regimes, such as high-speed flows and multiphysics turbulence, utilizing RANS and/or scale-resolving computational approaches.
Moreover, the strategy can be expanded by considering different SGD methods and optimization algorithms.
For instance, the Adadelta, Adam, SAG and SVRG SGD methods, and the Globally Convergent Method of Moving Asymptotes (GCMMA) optimization algorithm.
Finally, an avenue to further accelerate the optimization of aerodynamic systems by means of stochastic gradient-based approaches is to leverage the speedup of MF strategies based on combining the accuracy of high-fidelity models with the rapid computation of equivalent low-fidelity representations, as demonstrated in \cite{de2020bi}. 



\section*{Funding Sources}


This research was financially supported by the Beatriz Galindo Program (Distinguished Researcher, BGP18/00026) of the Ministerio de Ciencia, Innovaci\'on y Universidades, Spain. 

\section*{Acknowledgments}

AD acknowledges fruitful discussions with Prof. Kurt Maute (CU Boulder), Prof. John Evans (CU Boulder), and Dr. Subhayan De (CU Boulder) regarding this work. The authors additionally thank Dr. De for providing the MATLAB implementation of the AdaGrad algorithm.

\bibliography{References}

\end{document}